\documentclass[11pt,a4paper]{article}

\usepackage{amsthm, amsfonts,amssymb,euscript}
\usepackage{amsmath}
\usepackage{bm}
\usepackage{amssymb}
\usepackage{amsthm}
\usepackage{graphicx}
\usepackage{caption}
\usepackage{subcaption}
\usepackage{amsfonts}
\usepackage{float}
\usepackage{color}
\usepackage{slashed}
\usepackage[affil-it]{authblk}

\usepackage[usenames,dvipsnames,svgnames,table]{xcolor}
\usepackage[colorlinks=true]{hyperref}
\hypersetup{linkcolor=BrickRed, urlcolor=green, citecolor=blue, linktoc=page}

\numberwithin{equation}{section}

\newtheorem*{proposition*}{Proposition}
\newtheorem*{theorem*}{Theorem}
\newtheorem*{conjecture*}{Conjecture}
\newtheorem*{claim*}{Claim}
\newtheorem*{lemma*}{Lemma}
\newtheorem*{corollary*}{Corollary}
\newtheorem{theorem}{Theorem}[section]
\newtheorem{proposition}[theorem]{Proposition}

\newtheorem*{definition*}{Definition}
\newtheorem{definition}{Definition}[section]
\newtheorem*{assumption*}{\mathcal{A}ssumption}

\newtheorem*{remark*}{Remark}
\newtheorem{remark}{Remark}[section]

\numberwithin{equation}{section}
\begin{document}
\title{Instability results for transonic flows past airfoils}
\date{}
\author[1]{Yannis Angelopoulos \thanks {yannis@caltech.edu}}
\affil[1]{\small The Division of Physics, Mathematics and Astronomy, Caltech,
1200 E California Blvd, Pasadena CA 91125, USA}

\normalsize

\maketitle

\begin{abstract}
In this short note we present an instability result for transonic flows with respect to perturbations of the Mach number at infinity. More specifically we show that a perturbation of a transonic solution in the context of a Cauchy problem for the 2-dimensional steady, isentropic and irrotational Euler equations, which is a solution that solves an elliptic equation in certain parts of its domain and a hyperbolic equation in other parts, can only take place in the analytic category if the elliptic part of the domain is affected. 
\end{abstract}

\section{Introduction}
An interesting problem at the intersection of mathematics and aeronautics concerns the existence or not of special airfoils that allow flight at so-called transonic speeds -- that is speeds that are roughly between $0.8$ and $1.2$ times the speed of sound in vacuum -- without disturbances. 

It has been known to engineers that flight at ``slow" subsonic speeds is shock-free. At supersonic speeds on the other hand, the situation is quite different due to the presence of shocks. On a practical level, the presence or absence of shocks is important for the means of propulsion used for a flight, and of course for the cost required for it.

From a mathematical viewpoint, the problem is to study the \textbf{steady, isentropic and irrotational} Euler equations in a subset of $\mathbb{R}^2$. Modeling the case of flight at ``slow" subsonic speeds,  the problem reduces to the study of a second order quasilinear elliptic equation outside a convex obstacle (that is very similar to the minimal surface equation). This problem has been shown to be well-posed for ``generic" enough data, see the works \cite{bers-subsonic1}, \cite{bers-subsonic2}, \cite{shiffman} and the book of Bers \cite{bers-book} for a complete set of references. 

Modeling the case of supersonic flight in the same context results in the study of a second order hyperbolic Cauchy problem. In this case shocks are expected to form (in agreement with experimental results) but to our knowledge there are only heuristic arguments and numerical results towards this direction, see the book of Courant and Friedrichs \cite{courant-friedrichs} and the notes of Morawetz \cite{morawetz-notes}. A rigorous study of this situation will be addressed in upcoming work of the author.

Finally, the case of the so-called transonic flow is modeled  on an equation that is hyperbolic in a bounded subset of $\mathbb{R}^2$, and hyperbolic in the remaining part. Specific solutions have been constructed for this problem but the interesting problem from both a mathematical and engineering point of view is the study of stability of these solutions. This problem has been known as the \textit{transonic controversy} due to the differing opinions of physicists and engineers on the topic. For a very nice review of the topic see the important article of Morawetz \cite{morawetz-ams}, as well as the more recent \cite{morawetz-review-jhde}.

The controversy was settled to a large extent due to the work of Morawetz in \cite{morawetz1}, \cite{morawetz2}, \cite{morawetz3}, where she showed that the special transonic solutions that have been constructed so far are unstable with respect to ``generic" perturbations. She fully addressed the problem of perturbations in the hyperbolic part of the equation, but she only gave a partial result for perturbations at infinity (which is the problem that mostly interests the aerodynamicistics). Let us also mention that Morawetz's work is for symmetric profiles, but her results from \cite{morawetz1} have been extended to the non-symmetric case by Cook \cite{cook-transonic}.

\begin{figure}[H]
\begin{center}
\includegraphics[width=6cm]{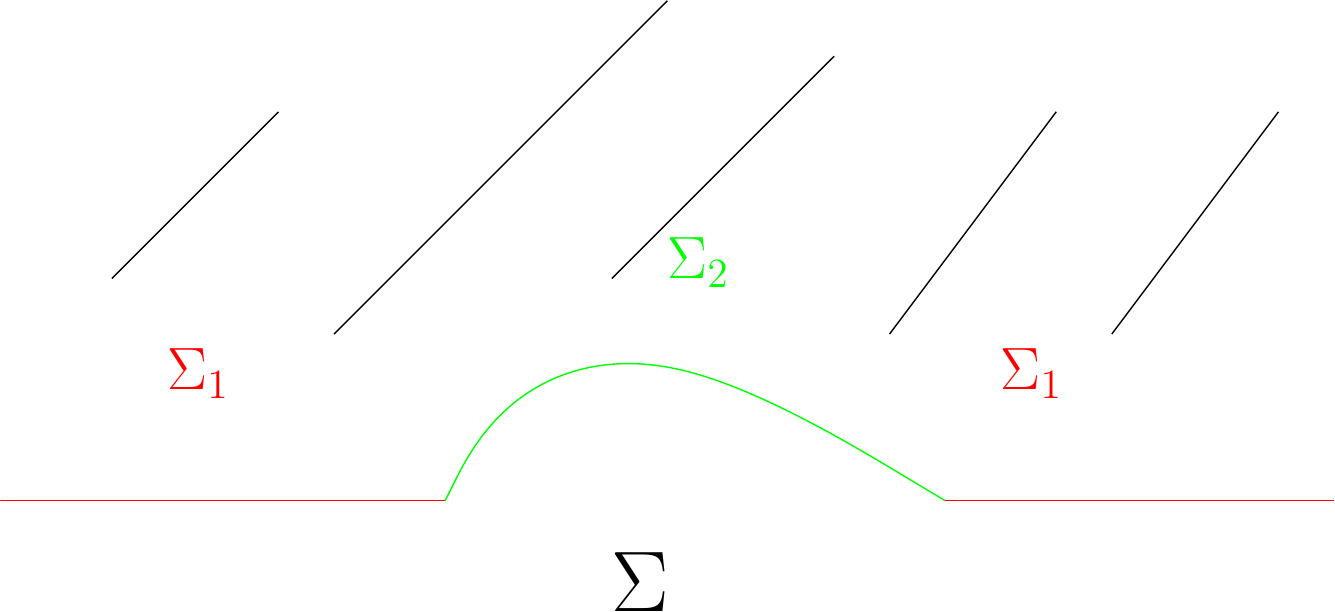}
\vspace{0.4cm}
\caption{\label{transonic-figure1}The initial hypersurface $\Sigma_{\tau}$, and the dotted area above it where we study our solutions.}
\end{center}
\end{figure}

In the current work we want to address the question of perturbation at infinity of transonic flows. The system of equations we are studying is a first order quasilinear system of the form:
$$ \partial_t u + A(u ) \cdot \partial_x u = 0 , $$
for $u : \mathbb{R}^2 \rightarrow \mathbb{R}^2$ and $A$ a $2\times 2$ matrix. We study a Cauchy problem where our initial hypersurface is a surface $\Sigma = \Sigma_1 \cup \Sigma_2$ as depicted in figure \ref{transonic-figure1}, and we solve the equation in the upper part of $\Sigma$ (i.e. the lined section in \ref{transonic-figure1}). Transonic solutions solve a hyperbolic problem in a neighbourhood of $\Sigma_2$, and an elliptic one in the remaining domain (for a more precise definition see Section \ref{transonic}).

\vspace{-0.8cm}

A rough statement of our main result is the following:
\begin{theorem}\label{thm:maininf}
Any perturbation of a transonic solution of the 2-dimensional steady, isentropic and irrotational Euler flow around a point of $\Sigma_1$ has to be a locally analytic perturbation.
\end{theorem}
Such perturbation include the case of perturbations of the Mach number at infinity, as such perturbations will inevitably include perturbations at the elliptic part of the initial data.

\subsection{Acknowledgements}
I would like to thank Toan Nguyen and Jared Speck for their interest in this work and for some helpful discussions about it.

\section{Overview}
In the next section we present the equation that we are working with. In Section \ref{transonic} we introduced the notion of transonic flows and review some past work on the topic. In Section \ref{Cauchy} we present several previous results on elliptic Cauchy problems. In Section \ref{proof} we state a precise version of our main Theorem along with its proof. Finally we include two appendices, Appendix \ref{appendixa} where we present some basic results on the FBI transformation, and Appendix \ref{matrix-pert} where we present some basic results on matrices with functional entries.

Let us note that the proof of the main Theorem of Section \ref{proof} is done in the spirit of the relevant Theorem in \cite{metivier-nonlinear}. More specifically we show the existence of a favorable representation for $BSu$, $u$ being our solution, where $B$ and $S$ are matrices, $B$ solving an elliptic boundary value problem, and $S$ the matrix whose columns are the eigenvectors corresponding to the complex eigenvalues of the matrix $A$. After showing that a certain set of points do not belong to the analytic wave front set of $BSu$ at the initial hypersurface $T=T_0$, we argue that the same should hold for $u$ by using now the aforementioned result.

\section{Compressible Euler equations in $(1+2)$-dimensions}

The isentropic compressible Euler equations on a 2-dimensional domain describe the motion of a perfect fluid with vanishing entropy. For $\rho$ the density, and $u = (u_1 , u_2 )$ the velocity vector we have the following set of equations:
\begin{equation}\label{eq:euler}
\begin{cases}
\frac{\partial \rho}{\partial t} + \frac{\partial ( \rho u_1 )}{\partial x_1} + \frac{\partial ( \rho u_2 )}{\partial x_2} = 0 \\
\frac{\partial u_i}{\partial t} + u_1 \frac{\partial u_i}{\partial x_1} + u_2 \frac{\partial u_i}{\partial x_2} = \frac{1}{\rho} \frac{\partial p}{\partial x_i}  \mbox{  for $i, j \in \{1,2\}$} ,
\end{cases}
\end{equation}
where $p \doteq p (\rho )$ is the pressure.

We will study the case of steady, isentropic and irrotational flow in a 2-dimensional domain. In the steady case the system \eqref{eq:euler} reduces to:
\begin{equation}\label{eq:euler_steady}
\begin{cases}
 \frac{\partial ( \rho u_1 )}{\partial x_1} + \frac{\partial ( \rho u_2 )}{\partial x_2} = 0 \\
 u_1 \frac{\partial u_i}{\partial x_1} + u_2 \frac{\partial u_i}{\partial x_2} = \frac{1}{\rho} \frac{\partial p}{\partial x_i}  \mbox{  for $i, j \in \{1,2\}$} \\
\frac{\partial u_1}{\partial x_2} = \frac{\partial u_2}{\partial x_1} , 
\end{cases}
\end{equation}
for functions 
$$ \rho: \mathcal{D} \rightarrow \mathbb{R}^{> 0} , \quad u : \mathcal{D} \rightarrow \mathbb{R} \times \mathbb{R} , $$
where $\mathcal{D} \subseteq \mathbb{R}^2$.

The assumption of the flow being irrotational implies that there exists a real valued function $\phi$ such that
$$ u \doteq \nabla \phi . $$
For a potential flow, Bernoulli's law implies that:
$$ \frac{u_1^2 + u_2^2}{2} + \int_1^{\rho} \frac{p' (s)}{s} \, ds = \mbox{  constant.  } $$
Differentiating the last equation in $\rho$ and using the first equation of \eqref{eq:euler_steady} we can eliminate the pressure and the density from the second equation of \eqref{eq:euler_steady}. From the last equation we also get that
$$ c^2 \doteq \frac{dp}{d\rho} = - \frac{\rho \cdot q}{\rho' (q)} , \quad M^2 = -\frac{q\cdot \rho' (q)}{\rho} \mbox{  for $q = \sqrt{u_1^2 + u_2^2}$,} $$ 
where by $\rho'$ we denote $\frac{d\rho}{dq}$, where $c$ denotes the \textit{local speed of sound}, and $M$ the \textit{Mach number} defined as
$$M \doteq \frac{q}{c} . $$

In the end we can reduce \eqref{eq:euler_steady} to the following first order quasilinear system of equations:
\begin{equation}\label{eq:1qe}
\frac{\partial}{\partial T} \begin{pmatrix} 
u_1 \\
u_2 
\end{pmatrix} + A(u ) \frac{\partial}{\partial x} \begin{pmatrix}
u_1 \\
u_2 
\end{pmatrix} = \begin{pmatrix}
0 \\
0
\end{pmatrix},
\end{equation}
setting $T = x_2$ and $x = x_1$, and where
\begin{equation*}
A(u) = \begin{pmatrix} 
- \frac{2u_1 u_2 }{c^2 - u_1^2} & \frac{c^2 - u_2^2}{c^2-u_1^2} \\
-1 & 0
\end{pmatrix} . 
\end{equation*}

A basic computation shows that the eigenvalues of $A$ are complex with non-zero imaginary part if
$$ q^2 = u_1^2 + u_2^2 < c^2 $$ 
and are real if
$$ q^2 = u_1^2 + u_2^2 > c^2 . $$

\section{Transonic flows past airfoils}\label{transonic}

One can study equation \eqref{eq:euler_steady} outside a convex obstacle $\mathcal{C}$ in $\mathbb{R}^2$, given also a boundary condition at infinity, this condition being the prescription of the Mach number at infinity, i.e.
$$ \lim_{|x| \rightarrow \infty} M \doteq M_{inf} . $$

For $M_{inf} < 1$ special transonic solutions have been constructed by various authors (see for instance the work of Tomotika and Tamada \cite{tomtam1}, \cite{tomtam2}, \cite{tomtam3}). Such solutions satisfy the following conditions:
\newline
1. $$  u \cdot  n_{\partial \mathcal{C}} = 0, $$
for $n_{\partial \mathcal{C}}$ the unit normal to $\partial \mathcal{C}$,
\newline
2. the flow is everywhere subsonic, i.e. $M < 1$, apart from two finite area regions, each one enclosed by  $\mathcal{P}' \subset \partial \mathcal{C}$ a connected component of the boundary of the convex obstacle $\mathcal{P}$, and a curve whose endpoints are the endpoints of $\mathcal{P}'$, where the flow is supersonic, i.e. $M >1$. Moreover, these two regions are symmetric with respect to the $x$-axis,
\newline
3. the flow is everywhere \textit{smooth} with the possible exception of a finite number of points on the boundary of the obstacle.

The question of interest (from an aerodynamics point of view) is the stability of such smooth solutions with respect to a variation of the airfoil which is taken to be the boundary of the convex obstacle $\mathcal{C}$ -- let us denote the airfoil by $\mathcal{P}  = \partial \mathcal{C}$. 

Morawetz in a series of works \cite{morawetz1}, \cite{morawetz2}, \cite{morawetz3} studied the stability of these solutions in the setting of an initial-boundary value problem. She studied equation \eqref{eq:euler_steady} with $u_1$, $u_2$ given on 
$$\Sigma \doteq \left(  \{ T = 0 \} \cap \{ | x | \geq c \} \right) \cup\mathcal{P}  , $$
and the solution is constructed in the upper half plane outside $\Sigma$. The reason that it is enough to study stability of the aforementioned smooth transonic flows in the upper half plane is by the fact that they are symmetric with respect to the $x$-axis. Morawetz showed that no \textbf{smooth} transonic flow $u'$ can exist that arises from data
\begin{equation}\label{pert1}
u_1' \cdot n_{\mathcal{P}_1} \neq 0 , \quad \left. u' \right|_{\Sigma \setminus \mathcal{P}_1}  = \left. u \right|_{\Sigma \setminus \mathcal{P}_1} , \quad \lim_{|x| \rightarrow \infty} M(u) = \lim_{|x| \rightarrow \infty} M(u') = M_{inf} < 1 , 
\end{equation}
for $\mathcal{P}_1 \subset \mathcal{P}'$. She also studied the case of possible smooth perturbations \textit{at infinity}, i.e.
$$ \lim_{|x| \rightarrow \infty} M(u' ) = M'_{inf} < 1 , $$
close to $M_{inf}$, which can be seen as a perturbation of the form
\begin{equation}\label{pert2}
u'_1 \cdot n_{\Sigma} = 0, \quad \left. | \partial_x^k ( u'_2 - u_2 ) | \right|_{\Sigma \cap \{ |x| \geq C\}} \leq \epsilon \mbox{  for all $k$,} \quad | M'_{inf} - M_{inf} | = \epsilon, \quad M'_{inf} < 1 , 
\end{equation}
for $\epsilon > 0$ small enough.

\begin{figure}[H]
\begin{center}
\includegraphics[width=9cm]{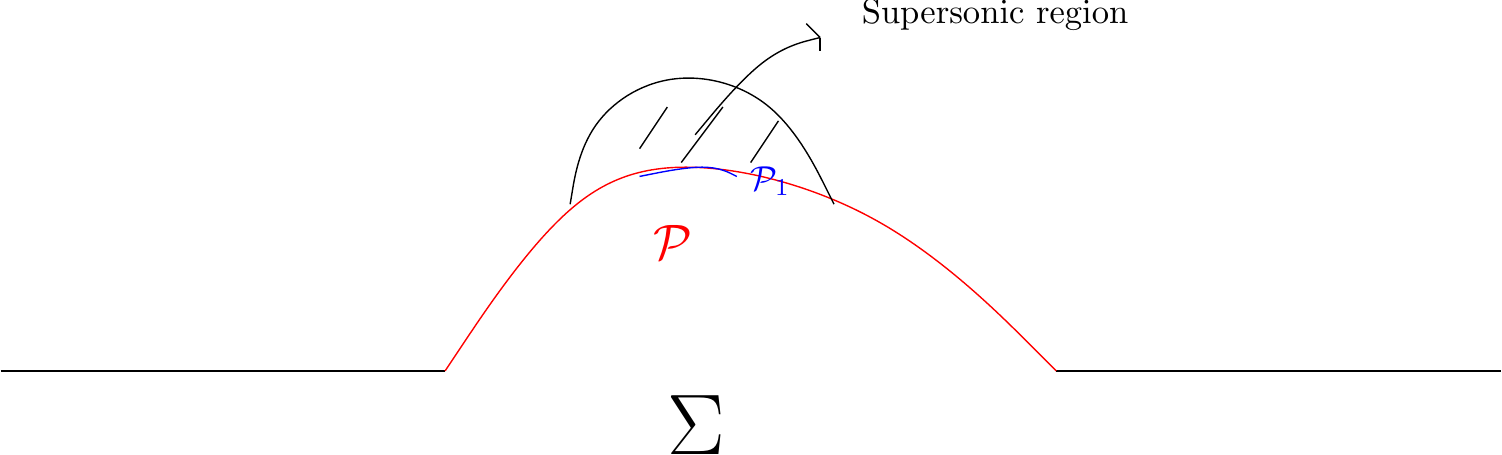}
\vspace{0.4cm}
\caption{\label{transonic-figure2}The initial hypersurface $\Sigma_{\tau}$, the profile $\mathcal{P}$, the subset $\mathcal{P}'$ of $\mathcal{P}$ where the perturbations considered by Morawetz in \cite{morawetz1} and \cite{morawetz2} are taken, and the supersonic bubble of a transonic solution.}
\end{center}
\end{figure}
\vspace{-0.8cm}

In all three of her aforementioned papers, Morawetz studied the steady, isentropic and irrotational Euler equations in 2D through the so-called hodograph transformation. This transformation can be considered as similar to the ones used in free boundary problems: in the present situation it is the sonic line that is fixed, i.e. the surface where $q = c$. It is a dynamical change of variables that reduces the system \eqref{eq:euler} to a second order quaslinear equation that is hyperbolic if $q > c$ and elliptic if $q < c$, which is studied in a subset of $\mathbb{R}^2$. 

In \cite{morawetz1} and \cite{morawetz2} Morawetz considers perturbations of type \eqref{pert1}, and in \cite{morawetz1} she shows instability for the linearized problem (see the figure \ref{transonic-figure2}), while in \cite{morawetz2} she shows instability for the full problem assuming that the nonlinear part remains small (which is a natural assumption). For perturbations of type \eqref{pert2} in \cite{morawetz3} she studied once again the linearized problem and she showed that if a solution to the perturbation problem exists, then there exists another transonic solution past a slightly airfoil $\mathcal{P}'$ with the Mach number at infinity of the perturbed solution obeying a slightly different density speed relation which is unstable with respect to perturbations at infinity. To our knowledge there is no result that shows instability directly for perturbation of type \eqref{pert2} for either the linearized or the full problem. 

All of Morawetz's results are proven through a contradiction argument using the vector field method (or what was called at the time the \textit{$abc$-method}). It is worth pointing out that the same argument was used earlier also by Frankl \cite{frankl1}.

Let us also note that Morawetz has attempted to construct transonic flows through viscocity methods, see the interesting paper of Gamba and Morawetz \cite{gamba-morawetz}. That program though to our knowledge remains incomplete. There are several other problems that someone can study as well, as the stability of transonic flows for the time-dependent problem. Note also that although transonic shocks have been constructed in other setting (see the work of Chen and Feldman \cite{chen-feldman}, and the works \cite{xin-yin1}, \cite{luo-xin}, \cite{luo-rauch-xie-xin}, as well as the references there-in).

In this note we want to address the non-existence of perturbations of special smooth transonic flows through a more general result that covers perturbation at the subsonic region of the initial hypersurface. Let us note at this point that the instabilities observed by Morawetz are due to the fact that the perturbations she considers make a mixed elliptic-hyperbolic problem overdetermined (for more examples on this see the book of Bers \cite{bers-book} and the one of Manwell \cite{manwell}), while the result that we show can be considered as a rigidity result.

\section{The nonlinear elliptic Cauchy problem}\label{Cauchy}
Considering perturbations at the subsonic region of the initial hypersurface, essentially reduces the problem to the question of well-posedness for a quasilinear locally elliptic Cauchy problem. There is a considerable amount of work (some of it being rather recent) on the well-posedness issue for elliptic Cauchy problems.

The first major result in the area is a celebrated theorem of Hadamard which states that the Cauchy problem for a \textit{linear second order elliptic equation} is ill-posed in the sense that for smooth non-analytic data there is no continuity for the data-to-solution map. It is worth mentioning that such a result is not true for analytic data due to the Cauchy-Kovalevskaya theorem, hence a corollary of the above is that any analytic solution of a linear second order elliptic Cauchy problem is unstable with respect to smooth perturbations.

On another direction, several results have been proven showing the necessity of hyperbolicity in various forms, for the well-posedness of a Cauchy problem, see \cite{metivier-nonlinear} and the references therein.

In the spirit of Hadamard's result, M\'{e}tivier in \cite{metivier-nonlinear}, showed that for first order quasilinear scalar equations of the form:
\begin{equation}\label{eq:metivier}
\begin{cases}
 \partial_t \psi +  \sum_{i=1}^d a_i (t,x,\psi) \partial_{x_i} \psi + b(t,x,\psi) = 0 , \\
  \psi (t_0 , x)  = f(x) ,
  \end{cases}
  \end{equation} 
where the $a_i$'s, $i \in \{1, \dots , d\}$ and $b$ are holomorphic functions of the arguments $t$, $x$ and $\psi : \mathbb{R} \times \mathbb{R}^d \rightarrow \mathbb{C}$ in a spacetime neighbourhood of $(t_0 , x , \psi (t_0 , x_0 ))$ for $t \geq t_0$, the Cauchy problem is ill-posed in the following ways: 

a) if equation \eqref{eq:metivier} has a $C^1$ solution in a spacetime neighbourhood of $(t_0 , x_0 )$ then $(x_0 , \xi ) \notin \textrm{WF}_A (f)$ for all $\xi \in \mathbb{R}^d$ such that $\sum_{i=1}^d \xi_i \Im a_i ( t_0 , x_0 , \psi (t_0 , x_0 )) > 0$,

b) the data-to-solution map in \eqref{eq:metivier} is not continuous with respect to the topology of Sobolev spaces.

The result of b) also holds for systems, i.e. for $\psi : \mathbb{R} \times \mathbb{R}^d \rightarrow \mathbb{C}^k$ for some $k \geq 1$ in \eqref{eq:metivier}, assuming that the matrices $a_i$, $i \in \{1,\dots , d\}$ have non-real eigenvalues close at $(t_0 , x_0 , \psi (t_0 , x_0 ))$, but the result of a) holds only for scalar equations and some very special cases of first order systems (essentially the ones that can be reduced to a scalar equation). The result of a) is of special interest to the current work. The method of proof involves a representation formula for the FBI transformation of the data $f$ (see Appendix \ref{appendixa} for the relevant definitions).

Following M\'{e}tivier's work, and using a variation of the method he introduced in order to show the result of a), Lerner, Morimoto and Xu in \cite{lerner-morimoto-xu} were able to relax the condition $\sum_{i=1}^d \xi_i \Im a_i ( t_0 , x_0 , \psi (t_0 , x_0 )) > 0$. Subsequently, Lerner, Nguyen and Texier in \cite{lerner-nguyen-texier} relaxed the conditions on the spectrum of the $a_i$'s, $i \in \{1, \dots , d\}$, for first order quasilinear systems and proved a result similar to the one of b). In that direction see also the recent works of Ndoumajoud and Texier \cite{ndoumajoud-texier1}, \cite{ndoumajoud-texier2}. 

\section{Non-existence of generic transonic solutions past airfoils for the initial-boundary value problem under smooth perturbations}\label{proof}

For $\Sigma$ given as in Section \ref{transonic}, being smooth, we will show the following result.

\begin{theorem}\label{thm:main1}
Assume that for some $(T_0 , x_0 ) \in \Sigma$ we have that for given \textit{analytic} initial-boundary data 
$$\left. (u_1 , u_2 )\right|_{\Sigma} , \quad \lim_{|x| \rightarrow \infty} M (u ) = M_{inf} ,$$ 
where
$$ A( \left. u \right|_{\Sigma} ) ( T_0 , x_0 ) \mbox{  has complex eigenvalues with non-zero imaginary part,} $$
the solution $u = (u_1 , u_2 )$ of equation \eqref{eq:euler_steady} that emanates from such data is \textbf{unstable} with respect to $C^{\infty}$ perturbations.

In particular, for a $C^1$ solution $u$ of equation \eqref{eq:euler_steady} and assuming that $\lambda_{\pm}$ are the two complex eigenvalues of $A$ at $(T_0 , x_0$ where
$$ \Im ( \lambda_+ ) = - \Im ( \lambda_- ) \neq 0, $$
we have that
$$ ( T_0 , x_0 , \xi ) \notin \mathrm{WF}_A ( \left. u \right|_{\Sigma} )\mbox{  for all $\xi \neq 0$,} $$
i.e. for $(T_0 , x_0 ) = y_0 \in \Sigma$, we have that $( y_0 , \xi )$ does not belong to the \textit{analytic wave front set} of $\left. u \right|_{\Sigma}$ for all $\xi \neq 0$.

\end{theorem}

\begin{remark}
Note that by our assumptions the eigenvalues of $A( \left. u \right|_{\Sigma} ) ( T_0 , x_0 )$ that we denote by $\lambda_+, \lambda_-$ have the form
$$ \lambda_+ = \kappa_1 + i \kappa_2 , \quad \lambda_- = \kappa_1 - i \kappa_2 , \mbox{  $\kappa_1 , \kappa_2 \in \mathbb{R}$, $\kappa_2 \neq 0$,} $$
and since $A$ is a $2\times 2$ matrix with real entries, we have that the corresponding eigenvectors have non-zero imaginary parts.
\end{remark}
\begin{remark}
The existence of a solution emanating from analytic data is a consequence of the Cauchy-Kowalevskaya theorem.
\end{remark}

\begin{proof}
Let us assume that we have a $C^1$ solution $u$ as in the assumptions of the Theorem. By our assumptions, there exists matrices $S$ and $D$, such that
$$ S^{-1} A S = D , $$ 
where $D$ is diagonal in a spacetime neighbourhood $\mathcal{W}_{T_0 , x_0}$ of $( T_0 , x_0 )$ having the form
$$ D (t,x)  = \begin{bmatrix} 
\lambda_1 (t,x,u) & 0 \\
0 & \lambda_2 (t,x,u)
\end{bmatrix}, $$
where $\lambda_1 (T_0 , x_0 ) = \lambda_+$, $\lambda_2 (T_0 , x_0 ) = \lambda_-$, for $\lambda_{\pm}$ the eigenvalues of $A$ at $(T_0 , x_0 )$. 

Note that then $S$ has the following form:
$$ S (t,x) = \begin{bmatrix} 
\sigma_1 (t,x,u) & s_1 (t,x,u) \\
\sigma_2 (t,x,u) & s_2 (t,x,u)
\end{bmatrix}, $$
where
$$ \begin{pmatrix}
\sigma_1 (T_0 , x_0 , \left. u \right|_{T=T_0}) \\
\sigma_2 (T_0 , x_0 , \left. u \right|_{T=T_0}) 
\end{pmatrix} , \quad \quad \begin{pmatrix}
s_1 ( T_0 , x_0 , \left. u \right|_{T=T_0} ) \\
s_2 ( T_0 , x_0 , \left. u \right|_{T=T_0} )
\end{pmatrix} , $$
are eigenvectors of $A$ at $(T_0 , x_0 )$ with respective corresponding eigenvalues $\lambda_+$, $\lambda_-$. These two representations for $D$ and $S$ follow from Proposition 1 of \cite{texier-matrix}, for $(t,x)$ in a sufficiently small spacetime neighbourhood of $(T_0 , x_0)$.

Note that from \eqref{eq:1qe} we have that
\begin{equation}\label{eq:2qe}
\partial_T ( Su ) + D \cdot \partial_x (Su) = [ ( \partial_T S ) \cdot S^{-1} ] \cdot (Su) + [ D \cdot ( \partial_x S ) \cdot S^{-1} ] \cdot (Su) .
\end{equation}

Let us now define a vector function
$$ \mathcal{Z} (t,x) \doteq \begin{pmatrix}
\zeta_1 (t,x,u) \\
\zeta_2 (t,x,u)
\end{pmatrix} , $$
such that for $(t,x)$ in the neighbourhood of $(T_0 , x_0 )$ where $D$ is diagonal, it satisfies the equation
\begin{equation}\label{eq:zeta}
\partial_t \mathcal{Z} + D \cdot \partial_x \mathcal{Z} = 0 ,
\end{equation}
with initial data
$$ \zeta_1 ( T_0 , x , , \left. u \right|_{T=T_0} ) = \zeta_2 ( T_0 , x  , , \left. u \right|_{T=T_0}) = x . $$

Let us also consider a $2 \times 2$ matrix $B$ to be specified later. Consider the following vector:
$$ BSu = \begin{pmatrix}
( BSu)_1 \\
( BSu )_2
\end{pmatrix} . $$
For $\lambda > 0$, $z \in \mathbb{C}$, and $q$ a quadratic form, we have by using equation \eqref{eq:2qe} that
\begin{equation}\label{eq:3qe}
\begin{split}
\partial_T & \begin{pmatrix} 
( BS u)_1 e^{-\lambda q ( \zeta_1 - z)} \\
( BS u)_2 e^{-\lambda q ( \zeta_2 - z)}
\end{pmatrix}  +  \partial_x \begin{pmatrix} 
( BS u)_1 e^{-\lambda q ( \zeta_1 - z)} \\
( BS u)_2 e^{-\lambda q ( \zeta_2 - z)}
\end{pmatrix} \\ = & \begin{pmatrix}
( [ \partial_T B + \partial_x (B \cdot D ) + B \cdot \partial_T S \cdot S^{-1} + B \cdot D \cdot \partial_x S \cdot S^{-1} ] Su )_1 e^{-\lambda q ( \zeta_1 - z)} \\
( [ \partial_T B + \partial_x (B \cdot D ) + B \cdot \partial_T S \cdot S^{-1} + B \cdot D \cdot \partial_x S \cdot S^{-1} ] Su )_2 e^{-\lambda q ( \zeta_2 - z)} 
\end{pmatrix}.
\end{split}
\end{equation}
Let us now consider the matrix equation
\begin{equation}\label{eq:B}
\begin{cases}
\partial_T B + \partial_x (B \cdot D ) + B \cdot \partial_T S \cdot S^{-1} + B \cdot D \cdot \partial_x S \cdot S^{-1} = 0 , \\
\left. B\right|_{\partial \mathcal{N}_{T_0 , x_0}} = B_d ,
\end{cases}
\end{equation}
which leads to our choice of $B$ close to $(T_0 , x_0 )$, where $\mathcal{N}_{T_0 , x_0} \subseteq \mathcal{W}_{T_0,x_0}$ is a spacetime neighbourhood of $(T_0 , x_0)$ where $A$ has compex eingevalues with non-zero impaginary part. Note that such a neighbourhood exists by Theorem \ref{thm:matrix2}. We can consider the data for $B$ on the boundary of a (spacetime) neighbourhood of $( T_0 , x_0 )$, that is $B_d$, to be analytic. Then the aforementioned equation is solvable as an elliptic boundary value problem, since $S$ and $D$ can be determined in $\mathcal{N}_{T_0 , x_0}$ by $u$ in $\mathcal{N}_{T_0 , x_0}$.

With this choice of $B$, equation \eqref{eq:3qe} becomes
\begin{equation}\label{eq:main}
\partial_T \begin{pmatrix} 
( BS u)_1 e^{-\mu ( \zeta_1 - z)^2 } \\
( BS u)_2 e^{-\mu ( \zeta_2 - z)^2}
\end{pmatrix} +  \partial_x \begin{pmatrix} 
( BS u)_1 e^{-\mu ( \zeta_1 - z)^2} \\
( BS u)_2 e^{-\mu ( \zeta_2 - z)^2}
\end{pmatrix}  = \begin{pmatrix}
0 \\
0
\end{pmatrix} .
\end{equation}
Then for $\chi$ a smooth cut-off function compactly supported in $(x_0 - R , x_0 + R)$ for some $R>0$ small enough, we have for some $T_0 < T_1$ that
\begin{equation}\label{eq:bsu}
\begin{split}
\int_{T=T_0} ( BSu )_i e^{-\mu ( \zeta_i - z)^2} \chi (x) \, dx = & \int_{T=T_1} ( BSu )_i e^{-\mu ( \zeta_i - z)^2} \chi (x) \, dx \\  & -\int_{T_0}^{T_1} \int_{\mathbb{R}} ( BSu )_i e^{-\mu ( \zeta_i - z)^2} \chi' (x) \, dx dT \mbox{  for $i \in \{1,2\}$.} \\ \doteq &  I_i + II_i \mbox{  for $i \in \{1,2\}$.}
\end{split}
\end{equation} 
We examine the terms separately. We start with the terms of type $I$ which have the form:
$$ I_i = \int_{T=T_1} ( BSu )_i e^{-\lambda ( \zeta_i - z)^2} \chi (x) \, dx, \mbox{  $i \in \{1,2\}$.} $$ 
We note that by the definition of the $\zeta_i$'s through equation \eqref{eq:zeta} and the assumption on the spectrum of $A$ we have that
$$ \Im \zeta_i (t,x)= -(t-T_0 ) \Im \lambda_i (0,x) + (t-T_0 ) \mathrm{r} (t,x) \mbox{  for $i \in \{1,2\}$ and $(t,x) \in \mathcal{N}_{T_0.x_0}$,} $$
where $\mathrm{r}$ is a remainder term for which we have that
$$ \lim_{t \rightarrow T_0^+} \sup_{y \in B(x_0 , r_0)} | \mathrm{r} | (t,y) = 0 , $$
for $r_0$ such that $\{t\} \times B (x_0 , r_0 ) \subseteq \mathcal{N}_{T_0 , x_0}$ for all $t \in [T_0 , T_1]$, assuming that $R < \frac{r_0}{2}$ so that the support of $\chi$ is included in $B(x_0 , r_0 )$.

First we note the following basic identity for a quadratic form acting on $\mathbb{C}$:
$$ \Re ( ) z-w )^2 ) = |Re z - Re w |^2 - | \Im z - \Im w |^2 . $$
From the above and the previous expression for $\Im \zeta_i$ we note that
\begin{align*}
\Re ( ( \zeta_i (t,x) - z) )^2 ) & \geq  - | \Im \zeta_i (t,x) - \Im z |^2 = - | - (t-T_0 ) \Im \lambda_i (0,x) + (t-T_0 ) r (t,x) - \Im z |^2 \\ = &   - | - (t-T_0 ) \Im \lambda_i (T_0,x) + (t-T_0 ) \Im \lambda_i (T_0 , x_0 ) \\ & - (t-T_0 ) \Im \lambda_i (T_0 , x_0 ) +  (t-T_0 ) r (t,x) - \Im z |^2 \\ \geq & -C (t - T_0 )^2 |  -  \Im \lambda_i (T_0,x) +  \Im \lambda_i (T_0 , x_0 ) |^2 - C |  - (t-T_0 ) \Im \lambda_i (T_0 , x_0 ) - \Im z |^2 \\ & - C (t-T_0 )^2 | r |^2 (t,x)) \\ \geq & -C (t-T_0 )^2 |  -  \Im \lambda_i (T_0,x) +  \Im \lambda_i (T_0 , x_0 ) |^2 \\ & - C | (t-T_0 ) \Im \lambda_i (T_0 , x_0 ) + \Im z |^2 - C  \cdot t^2 |r|^2 (t,x),
\end{align*}
for some $C > 0$ that depends on the area of the spacetime neighbourhood $\mathcal{N}_{T_0 , x_0}$ and again for $i \in \{1,2\}$. For some $t_0 \in [ T_0 , T_0 + c k^2 ]$ where $k \leq \frac{r_0}{2}$ and for $c > 0$ small enough such that $T_0 + ck^2 \leq C' \min ( r_0 , T_1 )$ for $C' > 0$ to be specified later, we suppose that we have that
$$ | (t_0 - T_0  ) \Im \lambda_i (0, x_0 ) + \Im z | \leq c' t_0 , $$
for some $c' > 0$ and for $i \in \{1,2\}$ (which is something that can always be arranged so -- we will make a choice for such $c'$ later). Then we have that
\begin{align*}
\Re ( ( \zeta_i (t_0,x) - z) )^2 ) + ( \Im z)^2  \geq & ( \Im z )^2 - C t_0^2 | x - x_0 |^2 \\ & - C |  - (t_0 - T_0 ) \Im \lambda_i (T_0 , x_0 ) - \Im z |^2 - C( t_0 - T_0 )^2 | r |^2 (t_0,x) \\ \geq & C'' | t_0 \Im \lambda_i (T_0,x_0 ) |^2   - C (t_0 - T_0 )^2 | x - x_0 |^2 \\ & - C    |(t_0 - T_0 ) \Im \lambda_i (T_0 , x_0 ) + \Im z |^2 - C (t_0 - T_0 )^2 | r |^2 (t_0 ,x) \\ \geq & C'' (t_0 - T_0 )^2 | \Im \lambda_i (T_0 , x_0 ) |^2 - C (t_0 - T_0 )^2 | x - x_0 |^2 \\ &  - C''' (t_0 - T_0 )^2 -  C | r |^2 (t_0 , x) ,
\end{align*}
for some constants $C'' , C''' > 0$ (again depending on the area of the spacetime neighbourhood $\mathcal{N}_{T_0 , x_0}$). Then we also have that:
\begin{align*}
& \inf_{\{ t, x , z |t \in [T_0, T_0 + ck^2 ],  | x - x_0 | \leq k , |\Im z - (t_0 - T_0 ) \Im \lambda_i (T_0 , x_0 ) | \leq c'' t_0\}}  \left[ \Re ( ( \zeta_i (t_0,x) - z) )^2 ) + ( \Im z)^2 \right] \\ \geq & C'' (t_0 - T_0 )^2 ( \Im \lambda_i (T_0 , x_0 )^2 - C (t_0 - T_0 )^2 | x - x_0 |^2 \\ &  - C''' (t_0 - T_0 )^2 -  C  | r |^2 (t_0 , x) .
\end{align*}
Now we can choose $k$ and $c'$ such that
\begin{align*}
C'' (t_0 - T_0 )^2 ( \Im \lambda_i (T_0 , x_0 )^2 & - C (t_0 - T_0 )^2 | x - x_0 |^2   - C'''(t_0 - T_0 )^2 \\ & -  C  | r |^2 (t_0 , x) \\ \geq & C'' (t_0 - T_0 )^2 ( \Im \lambda_i (T_0 , x_0 )^2 - C'''' | r |^2 (t_0 , x) ,
\end{align*}
for some $C'''' > 0$, and as for $t_0$ close to $T_0$ we have that 
$$ r (t_0 , x ) = o (t_0 - T_0 ) , $$
and in the end we have that
\begin{align*}
 \inf_{\{ t, x , z |t \in [T_0, T_0 + ck^2 ],  | x - x_0 | \leq k , |\Im z - t_0 \Im \lambda_i (T_0 , x_0 ) | \leq c'' t_0\}} & \left[ \Re ( ( \zeta_i (t_0,x) - z)^2 ) ) + ( \Im z)^2 \right] \\ \geq & C'' t_0^2 ( \Im \lambda_i (T_0 , x_0 ) )^2 . 
 \end{align*}
Gathering together all the above estimates, we have that
$$ I_i \leq \| (BSu )_i (T_1 ) \|_{ L^1 ( B(x_0 , r_0 )} e^{-C'' \mu t_0^2 ( \Im \lambda_i (T_0 , x_0 ))^2} , $$
for $i \in \{1,2\}$.

We now turn to the terms of type $II$ which have the form:
$$ II_i =  -\int_{T_0}^{T_1} \int_{\mathbb{R}} ( BSu )_i e^{-\mu | \zeta_i - z|^2} \chi' (x) \, dx dT \mbox{  for $i \in \{1,2\}$.} $$
In this case we cannot fully take advantage of the expression relating the $\zeta_i$'s to the $\lambda_i$'s, but on the other hand we can use the restricted range of the $x$ variable due to the presense of $\chi'$. Note that $x_0 \notin supp (\chi' )$. Since by the assumption on $\zeta_i$'s we have that:
$$ \Re (( \zeta_i (T_0 , x) - x_0 )^2) \geq C \cdot | x - x_0 |^2 , $$
for $C$ as before. By the fact that $\nabla ( \Re ( q (\zeta_i (t,x ) - z ))$ is continuous for 
$$ (t,x,z) \in [ T_0 , T_1 ] \times \bar{B}_{\mathbb{R}} (x_0 , r_0 / 2 ) \times B_{\mathbb{C}} ( x_0 , r_0 /2 ) , $$
for both $i \in \{1.2\}$ and for $T_1$ sufficiently close to $T_0$, we have that:
\begin{align*}
\inf_{(t,x,z) \in [T_0 , T_1 ] \times supp (\chi' ) \cap   \bar{B}_{\mathbb{R}} (x_0 , r_0 / 2 ) \times B_{\mathbb{C}} ( x_0 , r_0 /2 ) } & \Re ( |\zeta_i (t, x) - z |^2 )\\ \geq & C \frac{R^2}{4} - \bar{c} [ (t- T_0 ) + |z-x_0 | ] ,
\end{align*}
for some $\bar{c} > 0$. By choosing $\bar{c}'$ sufficiently small, we then have that
\begin{align*}
 \inf_{(t,x,z) \in [T_0 , \min ( T_0 + \bar{c}' R , T_1) ] \times supp (\chi' ) \cap   \bar{B}_{\mathbb{R}} (x_0 , r_0 / 2 ) \times B_{\mathbb{C}} ( x_0 , \min ( \bar{c}' R^2 , r_0 /2) ) }  & \Re ( ( \zeta_i (t, x) - z )^2) \\ \geq &  \bar{c}''  R^2 , 
 \end{align*}
for some $\bar{c}'' > 0$ that depends on the area of the spacetime neighbourhood $\mathcal{N}_{T_0 , x_0}$. 

Since we can choose $T_1$ as close to $T_0$ as we want, hence also choose $T_1 - T_0 \leq \bar{c}' R^2$, it is enough to bound instead of $II_i$ for $i \in \{1,2\}$, the quantities:
$$ III_i (s, z , \mu ) = -\int_{T_0}^{T_0 + s} \int_{\mathbb{R}} ( BSu )_i e^{-\mu q ( \zeta_i - z)} \chi' (x) \, dx dT \mbox{  for $i \in \{1,2\}$.} $$
Gathering together all the above estimates we have that
\begin{align*}
& \sup_{s \in [T_0 , \min ( T_0 + \bar{c}' R^2 , T_1 )], z \in B_{\mathbb{C}} ( x_0 , \min ( \bar{c}' R^2 , r_0 /2) )}   | III_i | (s, z , \mu)   \\ \leq &  \bar{C} R^3 \cdot \sup_{s \in [ T_0 ,  \min ( T_0 + \bar{c}' R^2 , T_1 )] , \frac{R}{2} \leq |x-x_0 | \leq  R} | (BSu )_i \cdot \chi' | \cdot e^{-\mu  \bar{c}'' R^2} ,  
\end{align*}
for some $\bar{C} > 0$. The above estimate along with the one for the quantities $I_i$, $i  \in \{1,2\}$, show that
\begin{equation}\label{est:basic}
\begin{split}
| & T \left. (BSu)_i \right|_{T=T_0} |  (z, \mu )  e^{- \Im z} \\ \leq & \mathcal{C} e^{-\epsilon \mu} \mbox{  for $i \in \{1,2\}$, for all $z$ in a (complex) neighbourhood of $x_0 - i T_0\Im \lambda_i (T_0 , x_0 )$, } 
\end{split}
\end{equation}
for some $\mathcal{C} , \epsilon >0$, and for $T \left. (BSu)_i \right|_{T=T_0}$ the FBI transformation of $ \left. (BSu)_i \right|_{T=T_0}$, $i \in \{1,2\}$, as defined in Appendix \ref{appendixa}. 

Assuming that 
$$\lambda_1 (T_0 , x_0 ) > 0 , \quad \quad \lambda_2 (T_0 , x_0 ) < 0 , $$
by the conic nature of the analytic wave front set, we get that 
$$ (x_0 , \xi ) \notin \mathrm{WF}_A ( \left. (BSu )_1 \right|_{T=T_0} ) \mbox{  for all $\xi > 0$, and  } (x_0 , \sigma ) \notin \mathrm{WF}_A (\left. ( BSu )_2 \right|_{T=T_0} ) \mbox{  for all $\sigma < 0$. } $$
Note now that since we are free to choose $B$ on $T_0$, we can claim that
$$ (x_0 , \xi ) \notin \mathrm{WF}_A ( \left. (Su )_1 \right|_{T=T_0} ) \mbox{  and  } (x_0 , \xi ) \notin \mathrm{WF}_A (\left. ( Su )_2 \right|_{T=T_0} ) \mbox{  for all $\xi \neq 0$. } $$

In order to make the same claim for $\left. u_1 , u_2 \right|_{T= T_0}$ we need to use the equations for both $u_1$, $u_2$ and $(Su )_1$, $(Su )_2$. 

We start by noticing that the previous process also shows that
$$ (x_0 , \xi ) \notin \mathrm{WF}_A ( \left. (Su )_1 \right|_{T=T_1} ) \mbox{  and  } (x_0 , \xi ) \notin \mathrm{WF}_A (\left. ( Su )_2 \right|_{T=T_1} ) \mbox{  for all $\xi \neq 0$, } $$
as once again we are free to choose $B$ on $T = T_1$. 

Moreover, we note that by repeating the previous process we have that
$$ (x_0 , \xi ) \notin \mathrm{WF}_A ( \left. (BSu )_1 \right|_{T=T'} ) \mbox{  for all $\xi > 0$, and  } (x_0 , \sigma ) \notin \mathrm{WF}_A (\left. ( BSu )_2 \right|_{T=T'} ) \mbox{  for all $\sigma < 0$, } $$
for any $T' \in [T_0 , T_1]$, and again by an appropriate choice of $B$, we can have that
$$ (x_0 , \xi ) \notin \mathrm{WF}_A ( \left. (Su )_1 \right|_{T=T'} ) \mbox{  and  } (x_0 , \xi ) \notin \mathrm{WF}_A (\left. ( Su )_2 \right|_{T=T'} ) \mbox{  for all $\xi \neq 0$, } $$
for any $T' \in [ T_0 , T_1 ]$. In order to see this in more detail we use once again the fact that $T_1$ can be chosen very close to $T_0$, and the properties of the $\zeta_i$'s, $i \in \{1,2\}$. We can just consider the quantities:
$$ \Big| \int_{\mathbb{R}} (Su )_i e^{-\mu ( x - z)^2} \chi (x) \, dx \Big| \mbox{  for $i \in \{1,2\}$,} $$
and we are able to show the required estimate by noticing that the lower bound for quantity in the exponential term is the same as for $i \in \{1,2\}$ we have that:
$$ \Re ( (x-z)^2 ) = \Re ( ( x - \zeta_i (t,x) + \zeta_i (t,x) - z )^2 ) \geq - ( \Im \zeta_i )^2 (t,x) - ( \Im \zeta_i (t,x) - \Im z )^2 , $$
and by choosing $T_1$ sufficiently close to $T_0$ the term $ - ( \Im \zeta_i )^2 (t,x)$ can be made small enough (due to the definition of the $\zeta_i$'s) so that it does not affect the rest of the analysis.

Next, we have the following equation for $u = (u_1 , u_2 )$ (which we write in terms of $S$ and $D$ and not $A$):
$$ \partial_T u + \partial_x ( S^{-1} \cdot D \cdot Su ) = [ \partial_x ( S^{-1} \cdot D ) + S^{-1} \cdot D \cdot ( \partial_x S ) \cdot S^{-1} ] Su . $$
We can multiply the above equation with $e^{-\frac{\mu}{2} (x-z)^2}$ for $z$ in an appropriately small neighbourhood of $x_0$ in $\mathbb{C}$, and we have the following equation:
\begin{equation}\label{eq:uaux}
\begin{split}
\partial_T ( u e^{-\frac{\mu}{2} (x-z)^2} ) +& \partial_x ( S^{-1} \cdot D \cdot Su e^{-\frac{\mu}{2} (x-z)^2} ) \\ = &  [ \partial_x ( S^{-1} \cdot D ) + S^{-1} \cdot D \cdot ( \partial_x S ) \cdot S^{-1} - \mu (x-z) \cdot ( S^{-1} \cdot D ) ] Su e^{-\frac{\mu}{2} (x-z)^2} .
\end{split}
\end{equation}
Integrating the above we get that
\begin{align*}
& \int_{T=T_0}  u e^{-\mu ( x- z)^2} \chi (x) \, dx =  \int_{T=T_1} u e^{-\mu ( x- z)^2} \chi (x) \, dx \\  & +\int_{T_0}^{T_1} \int_{\mathbb{R}} (S^{-1} \cdot D \cdot Su ) e^{-\mu ( \zeta_i - z)^2} \chi' (x) \, dx dT \\ &  + \int_{T_0}^{T_1} \int_{\mathbb{R}}  [ \partial_x ( S^{-1} \cdot D ) + S^{-1} \cdot D \cdot ( \partial_x S ) \cdot S^{-1} - \mu (x-z) \cdot ( S^{-1} \cdot D ) ] Su e^{-\mu (x - z)^2} \chi (x) \, dx dT \\ \doteq & i+ii+iii .
\end{align*}
Term $i$ can be treated similarly to terms $I_i$, $i \in \{1,2\}$ in \eqref{eq:bsu}, while term $ii$ can be treated as terms $II_i$, $i \in \{1,2\}$ in \eqref{eq:bsu}. For the last term we rely on the fact that the point of interest lies in the analytic wave front set of $Su$. We have that:
\begin{equation*}
\begin{split}
\Big| &\int_{T_0}^{T_1}  \int_{\mathbb{R}}   [ \partial_x ( S^{-1} \cdot D ) + S^{-1} \cdot D \cdot ( \partial_x S ) \cdot S^{-1} - \mu (x-z) \cdot ( S^{-1} \cdot D ) ] Su e^{-\mu ( x- z)^2} \chi (x) \, dx dT \Big| \\ \leq &  \int_{T_0}^{T_1} \sup_{x \in supp (\chi) }| \partial_x ( S^{-1} \cdot D )|(T,x) + |S^{-1} \cdot D \cdot ( \partial_x S ) \cdot S^{-1}| (T,x) + \mu (x-z) \cdot | S^{-1} \cdot D | (T,x ) \, dT \\ &\times \sum_{i=1}^2 \Big| \int_{\mathbb{R}} (Su )_i e^{-\mu ( x - z)^2} \chi (x) \, dx \Big| .
\end{split}
\end{equation*}
In the expression above the $\sup$ over the matrix quantities is meant as a $\sup$ over the sum of all their entries. 

Note that as $(x_0 , \xi ) \notin \mathrm{WF}_A (\left. (Su )_i \right|_{T=T'} )$ for $\xi \neq 0$, $i \in \{1,2\}$, and for all $T' \in [ T_0 , T_1 ]$, we have that
$$ \sum_{i=1}^2 \Big| \int_{\mathbb{R}} (Su )_i e^{-\mu ( x - z)^2} \chi (x) \, dx \Big| \leq \bar{\mathcal{C}} e^{\frac{\mu}{2} [ (\Im z )^2 - \epsilon ]} ,  $$
for some $\bar{\mathcal{C}} , \epsilon > 0$, and by Theorems \ref{thm:matrix1} and \ref{thm:matrix2} and the regularity of $u$, we have that:
$$  \int_{T_0}^{T_1} \sup_{x \in supp (\chi) }| \partial_x ( S^{-1} \cdot D )|(T,x) + |S^{-1} \cdot D \cdot ( \partial_x S ) \cdot S^{-1}| (T,x) + \mu (x-z) \cdot | S^{-1} \cdot D | (T,x ) \, dT \bar{\mathcal{C}}' , $$
for some $\bar{\mathcal{C}}' > 0$ depending on the area of the spacetime region $\mathcal{N}_{T_0 , x_0}$ and the $L^{\infty}$ norm of $u$ in $\mathcal{N}_{T_0 , x_0}$. Putting together the last two estimates we get the desired estimate for $iii$. Gathering together now all the estimates for $i$, $ii$, $iii$ we get that
$$ (x_0 , \xi ) \notin \mathrm{WF}_A (\left. u _i \right|_{T=T_0} ) \mbox{  for all $\xi \neq 0$}, $$
as desired.

\end{proof}

\section{Ill-posedness results in Sobolev and Gevrey spaces}
Let us note here that for initial data in Sobolev or Gevrey spaces we have the following result:
\begin{theorem}
For all $s \in \mathbb{R}$, $\alpha \in (0,1)$, $underline{r} > 0$, $\delta > 0$, there exists a sequence of positive numbers $\{ r_{\epsilon} \}$, converging to zero
$$ r_{\epsilon} \rightarrow 0 \mbox{  as $\epsilon \rightarrow 0$, } $$
such that for all initial data $\left. (u_1^{\epsilon} , u_2^{\epsilon} ) \right|_{\Sigma} \in H^s ( \Sigma) \times H^s ( \Sigma)$ for equation \eqref{eq:1qe} we have for the corresponding solutions $(u_1^{\epsilon} , u_2^{\epsilon} )$ that
$$ \lim_{\epsilon \rightarrow 0} \frac{\| u_1^{\epsilon} \|_{L^2 ( \Omega_{r_{\epsilon} ; \delta} )} + \| u_2^{\epsilon} \|_{L^2 ( \Omega_{r_{\epsilon} ; \delta} )}}{ ( \| \left. u_1^{\epsilon} \right|_{\Sigma} \|_{H^s ( B_{\underline{r}} (x_0 ) )} + \| \left. u_2^{\epsilon} \right|_{\Sigma} \|_{H^s ( B_{\underline{r}} (x_0 ) )} )^{\alpha} } = + \infty , $$
where $B_{\underline{r}} (x_0 )$ is the ball of radius $\underline{r}$ around $x_0$ in $\Sigma$ (where around $x_0$ at $\Sigma$ the matrix $A$ has complex eigenvalue with non-zero imaginary part), $H^s$ is the Sobolev space of order $s$, and where
$$ \Omega_{r_{\epsilon} ; \delta} \doteq \{ (t,x) \quad | \quad | x - x_0 |^2 + \delta t < r_{\epsilon}^2 \} . $$

Moreover, under the same assumptions as above for equation \eqref{eq:1qe} for all $\sigma \in ( 0 , 1/2 )$, $\alpha \in (0,1)$, $underline{r} > 0$, $\delta > 0$, there exists a sequence of positive numbers $\{ r_{\epsilon} \}$, converging to zero
$$ r_{\epsilon} \rightarrow 0 \mbox{  as $\epsilon \rightarrow 0$, } $$
such that for all initial data $\left. (u_1^{\epsilon} , u_2^{\epsilon} ) \right|_{\Sigma} \in G^{\sigma} ( \Sigma) \times G^{\sigma} ( \Sigma)$ for equation \eqref{eq:1qe} we have for the corresponding solutions $(u_1^{\epsilon} , u_2^{\epsilon} )$ that
$$ \lim_{\epsilon \rightarrow 0} \frac{\| u_1^{\epsilon} \|_{L^2 ( \Omega_{r_{\epsilon} ; \delta} )} + \| u_2^{\epsilon} \|_{L^2 ( \Omega_{r_{\epsilon} ; \delta} )}}{ ( \| \left. u_1^{\epsilon} \right|_{\Sigma} \|_{G^{\sigma ; c} ( B_{\underline{r}} (x_0 ) )} + \| \left. u_2^{\epsilon} \right|_{\Sigma} \|_{G^{\sigma ; c} ( B_{\underline{r}} (x_0 ) )} )^{\alpha} } = + \infty , $$
for $\Omega_{r_{\epsilon} ; \delta}$ as before, for $G^{\sigma}$ the Gevrey space of order $\sigma$, and for $\| \cdot \|_{G^{\sigma ; c} ( B_{\underline{r}} (x_0 ) )}$ the Gevrey norm given by
$$ \| f \|_{G^{\sigma ; c} ( B_{\underline{r}} (x_0 ) )} \doteq \sup_{\beta} \| \partial^{\beta} f \|_{L^{\infty} (B_{\underline{r}} (x_0 ) )} c^{-|\beta|}(  | \beta | ! )^{-1/\sigma} , $$
for functions $f$ such that
$$ \| \partial^{\beta} f \|_{L^{\infty} (B_{\underline{r}} (x_0 ) )} \leq C_{B_{\underline{r}} (x_0 ) } c_{B_{\underline{r}} (x_0 ) }^{|\beta|} (  | \beta | ! )^{1/\sigma} , $$
for $ C_{B_{\underline{r}} (x_0 ) }$, $c_{B_{\underline{r}} (x_0 ) } > 0$ depending on $B_{\underline{r}} (x_0 )$.
\end{theorem}

The aforementioned results shows that equation \eqref{eq:1qe} where the matrix $A$ is elliptic on a part of the initial hypersurface is ill-posed in Sobolev and Gevrey spaces in the sense that there is no continuity for the data-to-solution map. This result has been proven for general quasilinear first order systems. The Sobolev case was done by M\'{e}tivier \cite{metivier-nonlinear}, while the Gevrey case was proven by Morrise \cite{morisse}.
\appendix
\section{The analytic wave front set and some of its properties}\label{appendixa}
In this appendix we gather some facts on the analytic wave front set of a function. For a detailed exposition of the topics and a complete set of references we refer to the book of Delort \cite{delort-fbi}. First let us recall the following definition:
\begin{definition}
A point in phase space $(x_0 , \xi_0 ) \in T^* \mathbb{R}^d \setminus \{ 0 \}$ does \textbf{not} belong to the \textbf{analytic wave front set} of a distribution $f$, denoted by $WF_A (f)$, if there exists a neighbourhood $\mathcal{Z}$ of $z_0 \doteq x_0 - i \xi_0$ in $\mathbb{C}^d$ and some $\epsilon > 0$ such that
$$ \sup_{(z,\mu ) \in \mathcal{Z} \times [1,\infty)} e^{-\frac{\mu}{2} [ ( \Im (z) )^2 - \epsilon ]} | T f | ( z , \mu ) < \infty , $$
where $Tf$ is the \textbf{FBI transformation} of $f$ defined as:
$$ Tf (z , \mu) \doteq \int_{\mathbb{R}^d} e^{- \frac{\mu}{2} ( z-x )^2} f (x ) \, dx  \mbox{  for all $(z,\mu ) \in \mathbb{C}^d \times [ 1, \infty)$} . $$
Note that $Tf$ is a function on $\mathbb{C}^d \times [1, \infty)$ while $f$ is a function on $\mathbb{R}^d$.

\end{definition}
A basic change of variables (see \cite{lerner-morimoto-xu}) implies the following alternative characterization of the analytic wave front set:

\begin{proposition}
For $f$ a distribution, and $q$ a quadratic form on $\mathbb{C}^d$ given by
$$ q (z) \doteq \langle Qz , z \rangle \mbox{  for all $z \in \mathbb{C}^d$, } $$ 
for $Q$ a positive definite $d \times d$ matrix, we have that a point in phase space $(x_0 , \xi_0 ) \in T^* \mathbb{R}^d \setminus \{ 0 \}$ does \textbf{not} belong to $WF_A (f)$ the analytic wave front set of $f$ if and only if there exists a neighbourhood $\mathcal{N}_q$ of $z_0^q = x_0 - i Q^{-1} \xi_0$ in $\mathbb{C}^d$ and some $\epsilon > 0$ such that
$$ \sup_{( z , \mu ) \in \mathcal{N}_q \times [1,\infty)} e^{-\frac{\mu}{2} [ q ( \Im z ) - \epsilon ]} \left| \int_{\mathbb{R}^d} e^{-\frac{\mu}{2} q (z-x)} f(x) \,dx \right| < \infty . $$
\end{proposition}

\begin{remark}
Note that the quantity
$$ T_q f (z,\mu ) \doteq \int_{\mathbb{R}^d} e^{-\frac{\mu}{2} q (z-x)} f(x) \,dx ,$$
can be considered as a generalized FBI transformation relative to a quadratic form $q$. The aforementioned proposition implies that the analytic wave front set can be characterized in a the same way by the FBI transformation and by a generalized FBI transformation.
\end{remark}

In the case of functions, the FBI transformation (and hence its generalized form relative to a quadratic form $q$) can be used to characterize analyticity. We have the following Theorem for which a proof can be found in \cite{iagolnitzer} and in \cite{krantz-parks} for the case of $d=1$ that is of particular interest in the present work:
\begin{theorem}\label{thm:fbi-analytic}
A function $f$ is real analytic near a point $x_0 \in \mathbb{R}^d$ if and only if there exists neighbourhood $\mathcal{W}$ of $x_0$ and some $\Lambda > 0$ such that
$$ \left| \int_{\mathbb{R}^d} e^{-\frac{\mu}{2} ( x - y )^2 } e^{-2\mu \pi i x \cdot \xi} f(x) \, dx \right| \leq C e^{-\epsilon \mu} \mbox{  for all $y \in \mathcal{W},| \xi | > \Lambda$, $\mu > 0$,} $$
and for some positive constants $C$ and $\epsilon$. 

In particular we have that $f$ is real analytic near a point $x_0 \in \mathbb{R}$ if and only if 
$$(x_0 , \xi ) \notin \mathrm{WF}_A (f) \mbox{  for all $\xi \neq 0$.}  $$ 
\end{theorem}

\section{Some basic results of matrix perturbation theory}\label{matrix-pert}
We record here two basic theorems on matrices whose entries are functions. The first result relates the regularity of the entries of a matrix with the regularity of its eigenvalues and eigenvectors. For a proof see the note \cite{kazdan-note} of J. Kazdan.
\begin{theorem}\label{thm:matrix1}
Let $M(x)$ be a square matrix whose entries depend smoothly on $x$. For $\lambda$ a simple eigenvalue at $x=x_0$ with corresponding eigenvector $X_0$ we have that for all $x$ near $x_0$ there exists an simple eigenvalue $\lambda (x)$ of $M(x)$ with corresponding eigenvector $X(x)$ where both depend smoothly on $x$.

Moreover if the entries of $M(x)$ are real analytic functions of $x$, then the same conclusion as above holds with smoothness replaced by real analyticity, i.e. the eigenvectors and eigenvalues of $M(x)$ are real analytic in $x$.
\end{theorem}

The second result concerns the structure of the spectrum of a matrix with respect to a parameter. For a proof see \cite{texier-matrix}.
\begin{theorem}\label{thm:matrix2}
Let $M(x)$ be a square matrix whose entries are continuous with respect to $x$. Then the spectrum of $M$ is continuous with respect to $x$ in the following sense: for $\lambda_0 = \lambda (x_0 )$ an eigenvalue of $M(x_0)$ with multiplicity $m$, then there exists a neighbourhood $U$ of $x_0$ such that for all $x \in U$ the matrix $M(x)$ has $m$ eigenvalues $\lambda (x)$ (counting multiplicities) with $\lambda (x) \in B ( \lambda(x_0 ) , r )$ for some $r > 0$ (where $B(y,s)$ is the ball of radius $s$ centered at $y$).
\end{theorem}

\end{document}